\documentclass[10pt]{article}

\usepackage{amsmath}
\usepackage{amssymb}
\usepackage{amsfonts}
\usepackage{bm}
\usepackage{amsthm}
\newtheorem{theorem}{Theorem}

\newtheorem{definition}{Definition}
\newtheorem{lemma}{Lemma}

\theoremstyle{definition}
\newtheorem{example}{Example}

\usepackage{tikz}
\usepackage{pgfplots}
\usetikzlibrary{calc,plotmarks,patterns,decorations.text,intersections}
\usepackage{subcaption}
\usepackage{ifthen}

\newcommand{\alt}[1]{\ensuremath{\tilde{#1}}}		

\newcommand{\R}[1]{\ensuremath{\mathbb{R}^{#1}}}

\newcommand{\Zp}[1]{\ensuremath{\mathbb{Z}_{+}^{#1}}}

\newcommand{\B}[2]{\ensuremath{B_{#1}^{#2}}}
\newcommand{\pt}[2]{\ensuremath{\pi_{#1}^{#2}}}
\newcommand{\apt}[2]{\ensuremath{\alt{\pi}_{#1}^{#2}}}

\renewcommand{\d}{\ensuremath{d}}

\newcommand{\dt}{\ensuremath{\delta}}
\newcommand{\ord}{\ensuremath{r}}

\newcommand{\bset}[1]{\ensuremath{\mathbb{B}^{#1}}}
\newcommand{\bb}{\ensuremath{\boldsymbol{b}}}
\renewcommand{\b}{\ensuremath{b}}
\newcommand{\btset}[1]{\ensuremath{\mathbb{A}^{#1}}}
\newcommand{\bbt}{\ensuremath{\boldsymbol{a}}}
\newcommand{\bt}{\ensuremath{a}}

\renewcommand{\a}{\ensuremath{\alpha}}
\newcommand{\s}{\ensuremath{s}}
\renewcommand{\ss}{\ensuremath{\boldsymbol{s}}}

\renewcommand{\k}{\ensuremath{\kappa}}	
\newcommand{\kt}{\ensuremath{\lambda}}	
\newcommand{\ka}{\ensuremath{\sigma}}
\newcommand{\kb}{\ensuremath{\rho}}
\newcommand{\exam}{\ensuremath{\alpha}} 

\renewcommand{\l}{\ensuremath{\ell}}
\newcommand{\eps}{\ensuremath{\epsilon}}

\renewcommand{\u}{\ensuremath{\boldsymbol{u}}}
\renewcommand{\v}{\ensuremath{\boldsymbol{v}}}
\newcommand{\w}{\ensuremath{\boldsymbol{w}}}
\newcommand{\V}{\ensuremath{\mathcal{V}}}
\newcommand{\W}{\ensuremath{\mathcal{W}}}

\newcommand{\e}{\ensuremath{\boldsymbol{e}}}
\newcommand{\x}{\ensuremath{\boldsymbol{x}}}

\newcommand{\q}{\ensuremath{\boldsymbol{q}}}

\renewcommand{\c}{\ensuremath{c}}
\newcommand{\ct}{\ensuremath{c}}
\newcommand{\bnet}{\ensuremath{\mathcal{C}}}
\newcommand{\bnett}{\ensuremath{\mathcal{D}}}
\newcommand{\rais}[2]{\ensuremath{\mathrm{\boldsymbol{R}}_{#1}^{#2}}}
\newcommand{\lowe}[2]{\ensuremath{\mathrm{\boldsymbol{R}}_{#1}^{-#2}}}

\newcommand{\splo}[1]{\ensuremath{\Gamma^{#1}}}		
\newcommand{\simp}[2]{\ensuremath{\Delta_{#1}^{#2}}}		
\newcommand{\asplo}[1]{\ensuremath{\alt{\Gamma}^{#1}}}		
\newcommand{\asimp}[1]{\ensuremath{\alt{\Delta}^{#1}}}		
\renewcommand{\H}{\ensuremath{\mathrm{H}}}					

\newcommand{\zero}{\ensuremath{\boldsymbol{0}}}
\newcommand{\one}{\ensuremath{\boldsymbol{1}}}
\newcommand{\trans}[1]{\ensuremath{{#1}^{\intercal}}}

\newcommand{\aff}{\ensuremath{\mathrm{aff}}}
\newcommand{\conv}{\ensuremath{\mathrm{conv}}}

\newcommand{\fs}[1]{\scriptsize{#1}}

\renewcommand{\skip}{2.5pt}
\newcommand{\vertrad}{1.5pt}

\title{Towards the multivariate simplotope spline: continuity conditions in a class of mixed simplotopic grids}
\author{Tim~Visser\thanks{Department~of~Aerospace~Engineering, Delft~University~of~Technology, Delft, The~Netherlands} \and Cornelis~C.~de Visser \and Erik-Jan~van~Kampen}

\date{\today}

\begin{document}

\maketitle

\begin{abstract}


Smooth joins of simplex Bernstein-B\'{e}zier polynomials have been studied extensively in the past. In this paper a new method is proposed to define continuity conditions for tensor-product Bernstein polynomials on a class of mixed grids that meets certain out-of-facet parallelism criteria. The conditions are derived by first defining a simplex around the simplotopic bases of the tensor-product polynomials. Then the continuity conditions in the multivariate simplex spline defined on the resulting simplices, are adapted to hold for the tensor-product polynomials. The two- and three-dimensional results agree with the results found in the literature. It is expected that the method can be employed in more general grids.




\end{abstract}


\section{Introduction}
\label{sec:intr}









In this paper a new method is introduced for defining continuity conditions between tensor-product splines on $n$-dimensional simplotopes that share a complete facet. The latter implies the simplotopes are equal-dimensional and have at least one facet of equal shape. Although in this paper we describe a slightly more constrained case, the method is expected to be applicable to any pair of simplotopes that share a facet.

The research is sparked by innovations in the field of control theory. Many modern control techniques require high-dimensional function approximators, either for models of the system to be controlled \cite{Tol}, or for approximating the optimal control strategy \cite{Nithin}. In many cases a differentiable, smooth function is required. Over the last few years simplex splines have increasingly been used to fulfill this need \cite{CoenPhD}. Recently a more intuitive use of splines was proposed for finding the optimal control strategy in the reinforcement learning framework \cite{Nithin}.

In the past, splines of Bernstein-B\'{e}zier polynomials have been studied extensively. Although most research focused on the simplicial case \cite{Lai,ChuiBi}, some research has gone into splines on tensor-product bases and mixed grids \cite{LaiPhD,ChuiMu,Farin}. This research was however limited to the two- and three-dimensional cases. The contribution of this paper is an intuitive method for defining continuity conditions in mixed grids of any dimension. 

The geometric basis of tensor-product Bernstein polynomials is a product of simplices, often called a simplotope. Research into these polytopes is focused on finding triangulations for them \cite{DeLoera}, and its applications in game theory \cite{Freund,Talman}.

To combine the well-established theory of multivariate simplex splines and the general tensor-product basis, we redefine the simplotope as a specific subset of a higher-dimensional simplex. We show that if two simplotopes share a facet, so do the higher-dimensional simplices they define. Therefore we may construct continuity conditions in the higher dimensional simplex spline defined on the two simplices. These conditions are then adapted to hold for the lower-dimensional simplotope spline.

We start our discussion by thoroughly introducing the problem in section \ref{sec:prob}. Then in section \ref{sec:splo} the simplotope is introduced, which allows for a definition of the circumscribed simplex. In section \ref{sec:poly} the concept of the circumscribed simplex is extended to the polynomials defined on simplices and simplotopes. Finally in section \ref{sec:cont} the algorithm is presented for defining continuity conditions between tensor-product Bernstein-B\'{e}zier polynomials defined on a certain class of mixed grids of simplotopes. The results are summarized and discussed in section \ref{sec:conc}.








\section{Problem statement}
\label{sec:prob}

Let $\simp{}{n} \subset \R{n}$ denote an \emph{$n$-simplex}, defined as the convex hull of a vertex set $\V = \{ \v_{0},...,\v_n \}$. All basis polynomials will be expressed in barycentric coordinates with respect to $\V$, for which we introduce the set $\bset{n} = \{ \bb \in \R{n+1}: \sum \limits_{j=0}^{n} \b_j = 1 \}$. If a point $\x \in \R{n}$ has barycentric coordinates $\bb$ such that $\b_j \geq 0$ for all $j \in \{ 0,...,n\}$, then $\x \in \simp{}{n}$. If the equality holds for any $j$, $\x$ is on the boundary.

Let $\Zp{n}$ denote the set of $n$-element non-negative multi-integers. If $n=1$, the dimension will be omitted. We introduce for $\exam = (\exam_1,...,\exam_n) \in \Zp{n}$ the 1-norm $|\exam| = \exam_1 + ... + \exam_n$, the factorial $\exam! = \exam_1! \cdots \exam_n!$, and the inequality $\exam \geq 0$ to mean $\exam_j \geq 0, \forall j \in \{1,...,n\}$. Also, if we define $\x \in \R{n}$, then we write $\x^{\exam} = \prod_j \x_j^{\exam_j}$.

On the simplex $\simp{}{n}$ we may define a \emph{Bernstein basis polynomial} in terms of barycentric coordinates. For this we define the degree $\d \in \Zp{}$ and a multi-index $\k = (\k_0,...,\k_n) \in \Zp{n+1}$, where $|\k| = \d$. Then the basis polynomials are

\begin{equation}
\label{eq:prob_basi}
	\B{\k}{\d} (\bb) = \frac{\d!}{\k!} \bb^{\k} \text{, } \hspace{0.5cm} |\k| = \d.
\end{equation}

\noindent The complete set of basis polynomials for a given degree is defined by the valid permutations of $\k$. They are combined in a weighted sum to form a polynomial. The weights $\c_{\k}$, called \emph{B-coefficients}, have a spatial location in the simplex. They lie at the associated domain points $\q$ defined as

\begin{equation}
\label{eq:prob_bnet_simp}
\q = \sum \limits_{j=0}^{n} \frac{\k_{j}}{\d_j} \v_{j}.
\end{equation}

\noindent The collective of B-coefficients at their location in the simplex is called the \emph{B-net}.

If two simplices $\simp{}{n}$ and $\asimp{}{n}$ share $n$ vertices, they are said to share a facet. The single vertex that is not shared, is called the \emph{out-of-facet vertex}. Continuity can be enforced over the shared facet between Bernstein basis polynomials defined on the simplices. This is done by equating the directional derivatives of the polynomials on both sides of the shared facet in a direction $\u$ not parallel to this shared facet. Lai and Schumaker, among others, showed that this comes down to equating the De Casteljau iterations of the appropriate coefficients, with respect to directional coordinates $\ss$ describing the vector $\u$ \cite{Lai}. These conditions can then be manipulated further to come to the well known continuity conditions, but this is not done here. Then by assuming, without loss of generality, that the out-of-facet vertices are the first one in $\simp{}{n}$ and the last one in $\asimp{}{n}$, we find for continuity order $\ord$ the conditions

\begin{equation}
\label{eq:prob_cont}
	\c_{(0,\k_{1},...,\k_{n})}^{(\ord)} (\ss) =
		\alt{\c}_{({\k}_{1},...,{\k}_{n},0)}^{(\ord)} (\alt{\ss}), \hspace{0.5cm} 
		|\k| = |\alt{\k}| = \d-\ord,
\end{equation}

\noindent where $\c_{(0,\k_{1},...,\k_{n})}^{(\ord)} (\ss)$ signifies a De Casteljau iteration. Also, although the direction vector $\u$ is a single vector, the descriptions in directional coordinates $\ss$ and $\alt{\ss}$ may differ between the two simplices.

\begin{figure}
\centering
	\begin{tikzpicture}
		\input{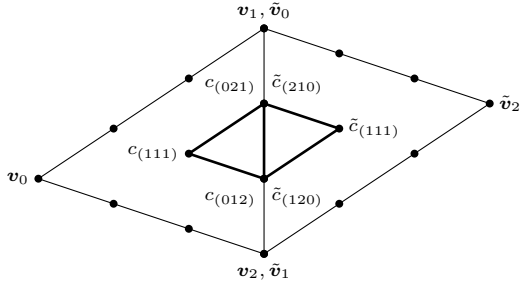}
	\end{tikzpicture}
	\caption{One of three first order continuity conditions in a cubic 2-simplex spline.}
	\label{fig:prob_simp}
\end{figure}

\begin{example}[Cubic spline on a triangulation]
\label{ex:prob_simp}
Consider two simplices $\simp{}{2}$ and $\asimp{}{2}$ that share an edge, as depicted in Figure \ref{fig:prob_simp}. The B-net of a cubic polynomial is indicated by the black circles in both simplices. In defining continuity conditions, we can choose any directional vector $\u$ not parallel to the shared edge $\{\v_1,\v_2\}$. If we choose $\u = \v_0 - \v_1$, the dashed first order continuity condition becomes

\begin{equation*}
\s_0 \c_{(111)} + \s_1 \c_{(021)} + \s_2 \c_{(012)} = \alt{\s}_{0} \alt{\c}_{(210)} + \alt{\s}_{1} \alt{\c}_{(120)} + \alt{\s}_{2} \alt{\c}_{(111)}.
\end{equation*}

\noindent This expression can be simplified by observing that $\u = (\v_0-\v_1) \parallel (\alt{\v}_1 - \alt{\v}_2)$. For the directional coordinates in $\simp{}{n}$ we have $\ss = (1,0,0)-(0,1,0) = (1,-1,0)$. In $\asimp{n}$ we have $\alt{\ss} = (1,1,-1)-(1,0,0)$, or due to the parallelism $\alt{\ss} = (0,1,0)-(0,0,1) = (0,1,-1)$. Filling this into the condition, we observe it simplifies to

\begin{equation*}
\c_{(111)} - \c_{(021)} = \alt{\c}_{(120)} - \alt{\c}_{(111)}.
\end{equation*}

\noindent The simplification lies in the fact that, most notably, $\alt{\v}_0$ is not required to describe the chosen vector $\u$. We will use this simplification in our final result. Note that a different vector choice leads to a different condition that is equally valid.
\end{example}

Now we combine simplices to form a new geometric basis.

\begin{definition}[Bottom-up simplotope]
\label{def:prob_bott}
Consider a multi-integer of dimensions $\nu = (\nu_1,...,\nu_\l) \in \Zp{\l}$. The \emph{$\nu$-simplotope} $\splo{\nu}$ is defined as the product of $\l$ simplices with dimensions $\nu$, that is $\splo{\nu} = \simp{1}{} \times ... \times \simp{\l}{}$, with $\simp{i}{}$ a simplex $\simp{}{\nu_i}$. 
\end{definition}

\noindent Note that the superscripted index in the simplex notation refers to the dimension of the simplex, whereas the subscripted one is simply an identifier. We will always assume the simplices are affinely independent, and therefore $\dim(\splo{\nu}) = |\nu|$. Also, we assume that the vertex sets $\W_i = \{\w_{ij}\}$ of the simplices $\simp{i}{}$ all contain the origin $\w_{i0} = \zero$. To simplify the notation, we define the vertices in these sets as $\w_{ij} \in \R{|\nu|}$, such that the product of simplices is not a concatenation, but a sum of the vertices of the separate simplices. The barycentric coordinates $\bbt_i \in \bset{\nu_i}$ of the parallel projections of a point onto the simplices $\simp{i}{}$, may be concatenated. Then we find a vector $\bbt \in \btset{\nu} = \{ \bbt \in \R{|\nu|+\l}: \sum \limits_{j=0}^{\nu_i} \bt_{ij} = 1, \forall i \in \{ 1,...,\l \} \}$. Note that we will generally have the index $i$ range from 1 to $\l$, and the index $j$ from 0 to $\nu_i$ or $n$. Finally, we say that two simplotopes $\splo{\nu}$ and $\asplo{\alt{\nu}}$ are of the same type if and only if $\nu = \alt{\nu}$.

\begin{example}[Simplotopes]
\label{ex:prob_bott}
In two dimensions there are two simplotopes, namely the parallelogram (or (1,1)-simplotope) and the 2-simplex (or (2,0)-simplotope). Note that we include the simplices in our definition, as a product of the $n$-simplex and enough zero-dimensional simplices to make $\l$ equal for both simplotopes. Although the product with simplices $\simp{}{0}$ does not affect the shape of the polytope, it is essential for a proper definition of the circumscribed simplex in section \ref{sec:splo}. In three dimensions we find three options: $\splo{(3,0,0)}$, $\splo{(2,1,0)}$ and $\splo{(1,1,1)}$, or in words the tetrahedron, the triangular prism and the parallelepiped respectively.
\end{example}

If a Bernstein basis polynomial is defined on each simplex $\simp{i}{}$ separately, we find a tensor-product Bernstein basis polynomial on the simplotope. To simplify the notation, we collect the degrees in a multi-integer $\dt \in \Zp{\l}$ and concatenate the multi-indices to $\kt = (\kt_1,...,\kt_\l) \in \Zp{|\nu|+\l}$, with $ \kt_i \in \Zp{\nu_i+1}, i \in \{1,...,\l \}$. Then the basis polynomials on the simplotope are

\begin{equation}
\label{eq:prob_tens}
	\prod \limits_{i=1}^{\l} \B{\kt_i}{\dt_i} (\bbt_i) = \frac{\dt!}{\kt!} \bbt^{\kt} 
		 \text{, } \hspace{0.5cm}|\kt_i| = \dt_i, \forall i \in \{1,...,\l\}.\\
\end{equation}

\noindent The complete set of tensor-product basis polynomials is defined by all possible combinations of valid permutations of the subsets $\kt_i$.

Given two simplotopes $\splo{\nu}$ and $\asplo{\alt{\nu}}$, with $|\nu| = |\alt{\nu}|$, we define tensor-product polynomials on both of them of the form

\begin{equation}
\label{eq:prob_poly}
\begin{split}
\pt{\nu}{\dt} (\bbt) &= \sum \limits_{|\kt_i| = \dt_i, \forall i }  
	\ct_{\kt} \frac{\dt!}{\kt!} \bbt^{\kt} \text{, and}\\
\pt{\alt{\nu}}{{\dt}} (\bbt) &= \sum \limits_{|\kt_i| = {\dt}_i, \forall i }  
	\ct_{\kt} \frac{{\dt}!}{\kt!} \bbt^{\kt}.
\end{split}
\end{equation}

\noindent The B-coefficients $\ct_{\kt}$ have a spatial location in the simplotope, at concatenations of domain points in the simplices $\simp{i}{}$. If we define $\W_i = \{ \w_{i0},...,\w_{i\nu_i}\}$ as the vertex set of $\simp{i}{}$, then these domain point are

\begin{equation}
\label{eq:prob_bnet_splo}
\q_i = \sum \limits_{j=0}^{\nu_i} \frac{\kt_{ij}}{\dt_i} \w_{ij}.
\end{equation}

\noindent The collective of these locations is called the B-net $\bnett$.

Our aim is to present an algorithm with which continuity conditions can be defined between the polynomials defined on a pair of simplotopes of different type. For this we will first discuss the shared facet of simplotopes. Then we refer to an elegant link between the simplotope and a higher-dimensional simplex, which is extended to the polynomials. Finally the continuity conditions in this higher-dimensional simplex spline are converted to hold on the simplotope.



%
%
%
%

\section{Properties of the simplotope}
\label{sec:splo}

Our approach will be to take continuity conditions from a higher-dimensional simplex spline, and adapt them to our simplotopic problem. For this we first need to find the shared facet between simplotopes, over which continuity conditions may be defined. Then we lay the link with the higher-dimensional simplex.

For two polytopes to share a facet (or ($n-1$)-face), they must have a facet that is identical between the two of them. The proper alignment of the vertices of these facets will only be considered later. First we formulate the following necessary condition for two simplotopes to share a facet.

\begin{lemma}[Shared facet]
If two simplotopes $\splo{\nu}$ and $\splo{\alt{\nu}}$, $\nu,\alt{\nu} \in \Zp{\l}$, share a facet, then there exist $\eps,\alt{\eps} \in \Zp{\l}, |\eps| = |\alt{\eps}| = 1$ such that $\nu - \eps = \alt{\nu} - \alt{\eps}$.
\end{lemma}

\begin{proof}
From the definition $\splo{\nu} = \simp{1}{} \times ... \times \simp{\l}{}$ we find that any facet of $\splo{\nu}$ can be found by reducing the (non-zero) dimension of one of the simplices in the product by one. The result is a lower-dimensional simplotope $\splo{\nu-\eps}$, with $ \eps \in \Zp{\l}: |\eps| = 1$. If the second simplotope has a facet of the same type, the simplotopes can share this facet.
\end{proof}

We learn from the above proof that the shared facet between simplotopes is again a simplotope. The simplex from the definition of $\splo{\nu}$ in which the vertex is excluded, that is where $\eps_i = 1$, plays an important role in defining continuity conditions.

\begin{definition}[Out-of-facet simplex]
\label{def:splo_outo}
Consider two simplotopes $\splo{\nu} = \simp{1}{} \times ... \times \simp{\l}{}$ and $\asplo{\alt{\nu}} = \asimp{1}{} \times ... \times \asimp{\l}{}$ that share, without loss of generality, a facet $\splo{\nu-\eps} = \simp{}{\nu_1-1} \times \simp{2}{} \times ... \times \simp{\l}{} = \asimp{1}{} \times ... \times \asimp{\l-1}{} \times \simp{}{\alt{\nu}_{\l}-1}$, where $\simp{}{\nu_1 -1} \subset \simp{1}{}$ and $\simp{}{\nu_\l -1} \subset \asimp{\l}{}$. Then $\simp{1}{}$ and $\asimp{\l}{}$ are the \emph{out-of-facet (oof) simplices} of $\splo{\nu}$ and $\asplo{\alt{\nu}}$ respectively.
\end{definition}

\noindent Note that, given a simplotope, we are free to choose the defining simplices $\simp{i}{}$. In this paper we will always choose them such that they share as many vertices as possible with the neighboring simplotope.

With the out-of-facet simplex in mind, we can formulate a condition that will prove to allow for a great simplification in the formulation of continuity conditions.

\begin{definition}[Out-of-facet cospatiality]
\label{def:splo_cosp}
Consider two simplotopes $\splo{\nu}$ and $\splo{\alt{\nu}}$ sharing a facet with out-of-facet simplices $\simp{1}{}$ and $\asimp{\l}{}$, where $\nu_1 \leq \alt{\nu}_\l$. $\splo{\nu}$ and $\asplo{\alt{\nu}}$ are said to be \emph{out-of-facet cospatial} if $\dim (\aff(\simp{1}{}) \cap (\asimp{\l}{} \setminus \simp{1}{})) \geq 1$
\end{definition}

\noindent In other words, two simplotopes are out-of-facet cospatial if there exists a vector in $\asimp{\l}{}$ that also lies in the affine hull of $\simp{1}{}$, but not on the shared facet.

\begin{figure}[t]
\centering
	\begin{tikzpicture}
		\input{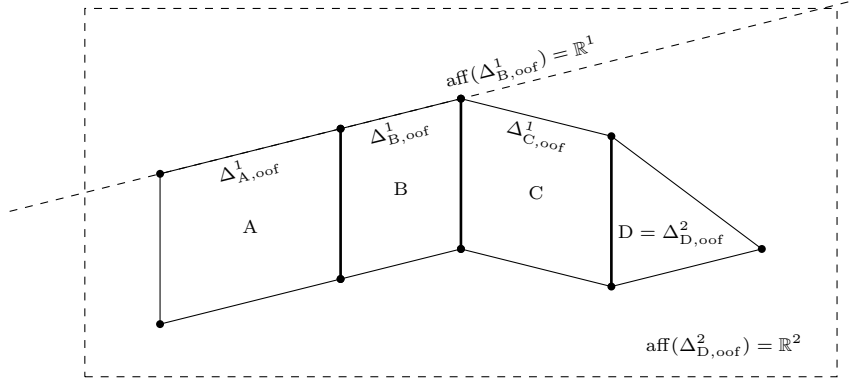}
	\end{tikzpicture}
	\caption{A tessellation in $\R{2}$ consisting of three parallelograms A, B, and C and a triangle D. The pairs \{A,B\} and \{C,D\} are out-of-facet cospatial, the pair \{B,C\} is not, for other pairs it is not defined.}
	\label{fig:splo_outo}
\end{figure}

\begin{example}[Out-of-edge simplices]
\label{ex:splo_outo}
Consider the tessellation in $\R{2}$ in Figure \ref{fig:splo_outo} consisting of three parallelograms A, B, and C and a triangle D. Clearly the pairs \{A,B\}, \{B,C\}, and \{C,D\} share 1-edges, indicated with thick lines. The out-of-facet (oof) simplices of all simplotopes are indicated. Note that $\simp{\mathrm{A,oof}}{1} \subset \aff(\simp{\mathrm{B,oof}}{1})$, so A and B are trivially out-of-facet cospatial. Clearly this does not hold for B and C, where $\aff(\simp{\mathrm{B,off}}{1}) \cap \simp{\mathrm{C,oof}}{1} \setminus \aff(\simp{\mathrm{B,off}}{1} \cap \simp{\mathrm{C,oof}}{1}) = \emptyset$. For the last pair, C and D, we note that $\simp{\mathrm{C,oof}}{1} \subset \aff(\simp{\mathrm{D,oof}}{2}) = \R{2}$. Therefore C and D are trivially out-of-facet cospatial. Any other pairs do not share a facet, and therefore their out-of-facet simplex is not defined.
%
\end{example}

Now refer back to the definition of barycentric coordinates $\bbt \in \btset{\nu}$ in the simplotope. Each subvector $\bbt_i$ describes the location in simplex $\simp{i}{}$, so that $|\bbt_i| = 1$. If we apply the conversion $\bb = \frac{1}{\l} \bbt \in \bset{|\nu|+\l-1}$, we find that these constraints describe hyperplanes with respect to a simplex $\simp{}{m}, m = |\nu|+\l-1$. This results in a definition of the simplotope as a subset of $\simp{}{m}$.

\begin{definition}[Top-down $\nu$-simplotope]
\label{def:circ_topd}
Consider a multi-index $\nu = (\nu_1,...,\nu_\l) \in \Zp{\l}$ and a simplex $\simp{}{m}, m=|\nu|+\l-1$ with vertices $\V$. Consider also a partition of $\V$ into $\l$ subsets $\V_i = \{ \v_{i0},...,\v_{i\nu_i} \}$ of $\nu_i+1$ vertices each. The \emph{$\nu$-simplotope} $\splo{\nu}$ is defined as the intersection $\simp{}{m} \cap \H_1 \cap ... \cap \H_\l$, with the hyperplanes $\H_{i} = \{ \bb \in \bset{|\nu|+\l-1}: \sum_j \b_{ij} = \a_i \}$ defined in barycentric coordinates with respect to $\V_i$, where $0 < \a_i < 1$ and $\sum_i \a_i = 1$. 
%
\end{definition}

\noindent Note that one of the hyperplanes is redundant, because inside $\simp{}{m}$ we have $|\bb| = 1$. One may therefore use $\l-1$ hyperplanes such that $\sum \a_i < 1$.

\begin{figure}[t]
\centering
	\begin{tikzpicture}
		\input{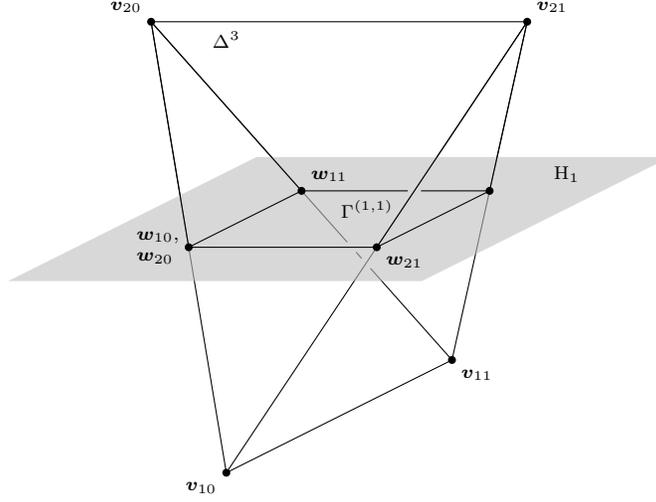}
	\end{tikzpicture}
	\caption{The parallelogram $\splo{(1,1)}$ is the intersection of the hyperplane $\H_1$ and the simplex $\simp{}{3}$.}
	\label{fig:circ_topd}
\end{figure}

\begin{example}[Parallelogram in a simplex]
\label{ex:circ_topd}
We will focus on the (1,1)-simplotope. From definition \ref{def:circ_topd} we may conclude that this polytope can be found in the ($|\nu|+\l-1 = 3$)-simplex $\simp{}{3}$. Describing the vertex set of $\simp{}{3}$ as $\V = \{ \v_{10},\v_{11},\v_{20},\v_{21}\}$, we can identify the two subsets $\V_{1} = \{ \v_{10},\v_{11}\}$ and $\V_2 = \{\v_{20},\v_{21}\}$ (see Figure \ref{fig:circ_topd}). Now the hyperplane $\H_1 =\{ \bb \in \bset{3}: \b_{10}+\b_{11} = \frac{1}{2} \}$ cuts the simplex in half, with $\V_1$ on one side and $\V_2$ on the other. The intersection $\H_1 \cap \simp{}{3}$ is the parallelogram $\splo{(1,1)}$. Note that $\H_2=\{ \bb \in \bset{3}: \b_{20} + \b_{21} = \frac{1}{2} \}$ trivially intersects the simplex in the same plane, as $\sum_{i,j} \b_{ij} = 1$.
\end{example}

Although proving the equivalence of the two definitions \ref{def:prob_bott} and \ref{def:circ_topd} is rather straightforward, we will limit ourselves to the part of the proof we need. 

\begin{lemma}[Circumscribed simplex]
\label{lm:circ_circ}
Given a simplotope $\splo{\nu}$ from definition \ref{def:prob_bott}, there exists a \emph{circumscribed simplex} $\simp{}{m}$ such that it defines $\splo{\nu}$ through definition \ref{def:circ_topd} with $\a_i = \frac{1}{\l}, \forall i$.
\end{lemma}

\begin{proof}
A simple construction of what we may call the \emph{standard circumscribed simplex} will suffice. Our approach is to drag certain subsets of vertices of $\splo{\nu}$ into the $(\l-1)$ extra dimensions in which $\simp{}{m}$ is defined. In this case we will assign each subset it's own extra dimension, except the last, which is simultaneously dragged into all extra dimensions.

As stated in section \ref{sec:prob}, we assume that each vertex set $\W_i$ describing the simplex $\simp{i}{}$ contains the origin $\w_0 = \zero \in \R{|\nu|}$. Then for the sets $\W_i, i \in \{ 1,...,\l-1 \}$ we scale the elements $\w_{ij} \in \R{|\nu|}$ with $\l$ and append the unit vector $\e_i \in \R{\l-1}$, to find $\v_{ij} = \trans{\left[ \l \trans{\w_{ij}},\trans{\e_{i}} \right]} \in \R{m}$. For $\W_\l$, we instead append the elements with $-\one \in \R{\l-1}$ (with $\one$ a vector of all ones) to find $\v_{\l j} = \trans{\left[ \l \trans{\w_{\l j}},\trans{-\one} \right]} \in \R{m}$. The complete set $\V = \{ \V_1 ,...,\V_\l\},\V_i = \{ \v_{i0},...,\v_{i\nu_i} \}$ is affinely independent, such that $\conv(\V) = \simp{}{m}$.

Now observe that the appended elements will cancel on the intersection of the hyperplanes $\H_i$ from definition \ref{def:circ_topd}, as $\frac{1}{\l} (-\one + \sum_i \e_i) = \zero$. What remains on the intersection is a linear combination of the vertices $\l \w_{ij} \in \R{|\nu|}$, with weights $\bb \in \frac{1}{\l} \btset{\nu}$. The factors of $\l$ cancel, and we are left with the vertices $\w_{ij}$ and weights $\bbt \in \btset{\nu}$.
\end{proof}

Note that the vertex sets $\V_i$ are parallel, scaled copies of the sets $\W_i$. By translating these vertex sets appropriately, infinitely many circumscribed simplices can be constructed, given a simplotope.

\begin{example}[Circumscribed tetrahedron]
\label{ex:circ_circ}
We will construct the simplex $\simp{}{3}$ from example \ref{ex:circ_topd} using the standard circumscription procedure used above. This time we start from the (1,1)-simplotope $\splo{(1,1)}$. It is defined by the product $\simp{1}{1} \times \simp{2}{1}$ with vertex sets $\W_1 = \{ 0,1 \}$ and $\W_2 = \{ 0,1 \}$. That is, $\splo{(1,1)}$ has vertices $\{ \}$. Appending $\e_1 \in \R{}$ and $-\one \in \R{}$ to these sets respectively, and scaling with $\l$, we find $\V_1 = \{ \v_{10},\v_{11} \} = \{ \trans{(0,0,1)},\trans{(2,0,1)}\}$ and $\V_2 = \{ \v_{20},\v_{21} \} = \{\trans{(0,0,-1)},\trans{(0,2,-1)} \}$. Clearly $\conv(\V) = \simp{}{3}$ and $\{ \H_1: \b_{10} + \b_{11} = \frac{1}{2} \} \cap \simp{}{3} = \splo{(1,1)}$.
\end{example}

In order to deduce the continuity conditions between simplotopes from those of the circumscribed simplices, the simplices need to share a facet. Therefore we prove that if two simplotopes share a facet, we can construct simplices that do too.

\begin{lemma}[Pair of circumscribed simplices]
\label{lm:circ_shar}
Consider two simplotopes $\splo{\nu}$ and $\asplo{\alt{\nu}}$ that share a facet $\splo{\nu-\eps}, \eps \in \Zp{\l}, |\eps| = 1$. There exist circumscribed simplices $\simp{}{m}$ and $\asimp{}{m}$ for $\splo{\nu}$ and $\asplo{\alt{\nu}}$ respectively, such that $\simp{}{m} \cap \asimp{}{m} = \simp{}{m-1}$, with $\simp{}{m-1}$ a facet of both simplices.
\end{lemma}

\begin{proof}
$\splo{\nu}$ and $\asplo{\alt{\nu}}$ share all but one simplex completely, the exception being the out-of-facet simplex. Without loss of generality we may assume that $\w_0 = \zero$ is a vertex of all of the shared simplices, and that $\w_{1\nu_1}$ and $\w_{\l \nu_\l}$ are the out-of-facet vertices. Then the standard circumscribed simplex of both simplotopes will include the vertices $\v_{ij} = \trans{\left[ \l \trans{\w_{ij}},\trans{\e_{i}} \right]}$ for $i \in \{2,...,\l-1\}, j \in \{0,...,\nu_i\}$ and $i = 1, j \in \{0,...,\nu_1-1\}$, and $\v_{\l j} = \trans{\left[ \l \trans{\w_{\l j}},\trans{\zero} \right]}$ for $j \in \{ 0,...,\nu_\l \}$. In total the circumscribed simplices will therefore share $|\nu-\eps| + \l - 1 = m-1$ vertices.
\end{proof}

As before, the circumscribed simplices are not unique. Pure translation of the vertex sets will provide an infinite amount of solutions. When the circumscribed simplex of the shared facet is set however, the remaining two out-of-facet vertices are uniquely defined.

\section{Equivalence of simplex and simplotope polynomials}
\label{sec:poly}

We have shown in the previous section that simplotopes are strongly connected to a higher-dimensional simplex. Now we will do the same for the polynomials defined on both polytopes.

Consider a tensor-product basis polynomial as in \eqref{eq:prob_tens}. The barycentric coordinates can be replaced by their simplex counterparts by filling in $\bbt = \l \bb$, to find

\begin{equation}
\label{eq:poly_equi}
\begin{split}
	\frac{\dt!}{\kt!} \bbt^{\kt} & = \frac{\dt!}{\kt!} \l^{|\dt|} \bb^{\kt} \\
		& = \frac{\dt!}{|\dt|!} \l^{|\dt|} \B{\kt}{|\dt|}(\bb).
\end{split}
\end{equation}

\noindent In words, the tensor-product polynomials are a scaled cut of Bernstein basis polynomials of a higher dimension (namely $m = |\nu|+\l-1$), defined by the constraints on the local barycentric coordinates $\sum_j \bt_{ij} = 1$.

\begin{example}[Bi-quadratic polynomial]
\label{ex:poly_biqu}
If we combine $\B{(2,0)}{2}$ and $\B{(1,1)}{2}$ we find the two-dimensional bi-quadratic polynomial

\begin{equation*}
	\B{(2,0)}{2} (\bbt_1) \B{(1,1)}{2} (\bbt_2) = 2 \bt_{10}^{2} \bt_{20} \bt_{21}.
\end{equation*}

\noindent By filling in $\bt_{ij} = \l \b_{ij}$ and $\l=2$, we may write the above as 

\begin{equation*}
	\B{(2,0)}{2} (\bbt_1) \B{(1,1)}{2} (\bbt_2) = 32 \b_{10}^{2} \b_{20} \b_{21}.
\end{equation*}

\noindent We may compare this to the ($|\nu|+\l-1 = 3$)-dimensional Bernstein basis polynomial of degree $|\dt| = 4$

\begin{equation*}
	\B{(2,0,1,1)}{4} (\bb) = 12 \b_{0}^{2} \b_{2} \b_{3},
\end{equation*}

\noindent to find that they differ a factor of $\frac{\dt!}{|\dt|!} \l^{|\dt|} = \frac{8}{3}$.
\end{example}

A similar link can be found between the B-nets of the simplex and tensor-product polynomials. From the above it is clear that a valid multi-index $\kt$ for a polynomial of degrees $\dt$, is also a valid multi-index for the polynomial of dimension $m = |\nu|+\l-1$ and degree $|\dt|$. Because each multi-index $\kt$ is coupled directly to a B-coefficient, we can find the B-net of the tensor-product polynomial in that of the $m$-dimensional simplex polynomial of total degree.

\begin{lemma}[B-net equivalence]
\label{lm:poly_bnet}
Consider a tensor-product polynomial of degrees $\dt$ with a B-net $\bnett$ defined on a simplotope $\splo{\nu}$ with circumscribed simplex $\simp{}{m}$. If a Bernstein-B\'{e}zier polynomial of total degree $\d = |\dt|$ with B-net $\bnet$ is defined on $\simp{}{m}$, then $\bnett = \bnet \cap \splo{\nu}$, with $\splo{\nu}$ defined through definition \ref{def:circ_topd} with the same vertex partition and with weights $\a_i = \frac{\dt_i}{|\dt|}$.
\end{lemma}

\begin{proof}
The B-coefficients $\c_{\k} \in \bnet$ lie at domain points $\q$ as defined in \eqref{eq:prob_bnet_simp}. In other words, they lie at locations with barycentric coordinates $\bb =\frac{1}{\d} \k$. Then if $\k = \kt$, the constraints from the tensor-product polynomial $|\kt_i| = \dt_i$ define subsets of the hyperplanes $\H_i = \{ \bb \in \bset{m}: \sum_j \b_{ij} = \frac{\dt_i}{\d} \}$. $\bnett$ consists of those B-coefficients for which all constraints are satisfied simultaneously, that is $\bnett = \bnet \cap \H_1 \cap ... \cap \H_\l $.
\end{proof}

\begin{example}[Tensor-product B-net]
\label{ex:poly_topd}
In Figure \ref{fig:cont_retr_bnet} the B-net of two simplex polynomials is displayed for degree $\d = 4$. According to lemma \ref{lm:poly_bnet}, the B-nets of degrees $\dt = (2,2)$ can be found in this B-net at the simplotope as in definition \ref{def:circ_topd} with weights $\a_1 = \a_2 = \frac{1}{2}$, defined with respect to the same vertex sets as depicted in Figure \ref{fig:cont_retr_base}. This B-net is shown in the left simplex in Figure \ref{fig:cont_retr_splo}. Similarly, the B-net of degrees $\dt = (1,3)$ can be found at the simplotope with weights $\a_1 = \frac{1}{4}$ and $\a_2 = \frac{3}{4}$. This B-net is depicted in the right simplex in Figure \ref{fig:cont_retr_splo}. Note that between the two simplices, the partition of the vertices into sets is different. In Figure \ref{fig:cont_retr_fina}, the B-net of a bi-quadratic polynomial is highlighted in both simplices.
\end{example}

\section{Continuity conditions between out-of-facet cospatial simplotopes}
\label{sec:cont}

Let $\splo{\nu}$ and $\asplo{\alt{\nu}}$ be two equal-dimensional simplotopes that share a facet and are out-of-facet cospatial. Assume that we have added zero-dimensional simplices, if necessary, so that the simplotopes are defined by the same amount of simplices $\l$. Without loss of generality we assume that the out-of-facet simplices are $\simp{1}{}$ and $\asimp{\l}{}$ respectively, and that $\nu_1 < \alt{\nu}_\l$. That is, defining a multi-index $\eps^{(i)}$ with $ \eps_{i}^{(i)} = 1$ the only non-zero element, we have $\nu - \eps^{(1)} = \alt{\nu} - \eps^{(\l)}$. Tensor-product polynomials $\pt{\nu}{\dt}$ and $\apt{\alt{\nu}}{\dt}$ with equal degrees are defined on these simplotopes. 

Before formulating our main result, we first introduce the degree raising operator $\rais{}{k}$ \cite{Lai}. It is used to write a coefficient as the convex combination of coefficients of a lower degree, equal-dimensional polynomial as

\begin{equation}
\label{eq:cont_rais}
	\rais{}{k} \c_{\k} = \frac{\d!}{(\d+k)!} \sum \limits_{|\mu|=k} \frac{\k!}{(\k-\mu)!} \c_{\k-\mu}, \hspace{0.5cm}
	|\k| = \d+k.
\end{equation}

\noindent Because generally we will raise the degree within a simplex of a simplotope, we introduce the local degree raising operator as

\begin{equation}
\label{eq:cont_lora}
	\rais{i}{k} \ct_{\kt} = \frac{\dt_i!}{(\dt_i+k)!} 
		\sum \limits_{|\mu|=k} \frac{\kt!}{(\kt-\mu)!} 
		\c_{(\kt_1,...,\kt_{i-1},\kt_i-\mu,\kt_{i+1},...,\kt_\l)}, \hspace{0.5cm}
	|\kt_i| = \dt_i + k.
\end{equation}

\noindent The inverse operator $\lowe{i}{k}$ will be used to describe coefficients of a lower degree polynomial by a linear combination of higher degree coefficients \cite{Farouki}. For general degree $\d$ and change of degree $k$ it can be found through iteration on \eqref{eq:cont_rais} to be

\begin{equation}
\label{eq:cont_lowe}
\begin{split}
	\lowe{}{k} \c_{\k} & = \sum \limits_{t=0}^{\d-\k_j}
		\sum \limits_{|\ka| = t}
		(-1)^t \frac{t!}{\ka!} \frac{(k-1+t)!}{(p-1)!t!}
		\frac{(\d+k)!}{(\k+\kb)!}
		\frac{\k!}{\d!}
		\c_{\k+\kb},
		\\
	\text{with } \kb & = (-\ka_0,...,-\ka_{j-1},t+k,-\ka_{j+1},...,-\ka_{n}), 
		\text{ and } \ka_{j}=0.
\end{split}
\end{equation}

\noindent Like degree raising, reduction can be done within a single simplex $\simp{i}{}$.

Now we can present our main theorem.

\begin{theorem}[Continuity in mixed grids]
\label{th:cont_cont}
An $\ord^{\text{th}}$ order smooth join between $\pt{\nu}{\dt}$ and $\apt{\alt{\nu}}{\dt}$ can be established by enforcing the conditions

\begin{equation}
\label{eq:cont_splo}
\left( \ct_{(\kt_{11},...,\kt_{1, \nu_1 -1},0)}^{(\ord)} (\ss_1) \right)_{(\kt_{2},...,\kt_{\l})} =
	\lowe{1}{\ord} \rais{\l}{\ord} \left( \alt{\ct}_{(\kt_{11},...,\kt_{1, \nu_1 -1},\kt_{2},...,\kt_{\l-1})} \right)_{(\kt_{\l 0},...,\kt_{\l \nu_\l},0)}^{(\ord)} (\alt{\ss}_{\l}),
\end{equation}

\noindent for all $\kt$ for which $|\kt_i| = \dt_i$ and $\kt_{1 \nu_1} = \kt_{\l \alt{\nu}_\l} = \ord$.
\end{theorem}


\begin{proof}
We start by defining the circumscribed simplices $\simp{}{m}$ and $\asimp{}{m}$ according to lemma \ref{lm:circ_shar} such that they share a facet. Note that definition \ref{def:circ_topd} implies that the vertex subsets $\V_i$ from the proof of lemma \ref{lm:circ_circ} are parallel to the simplices $\simp{i}{}$. In combination with the out-of-facet cospatiality of $\splo{\nu}$ and $\asplo{\alt{\nu}}$ this means we may choose the vector $\u$ such that it is parallel to both $\V_1$ and $\alt{\V}_{\l}$. If we assume, without loss of generality, that $\u = \v_{1,\nu_1} - \v_{1,\nu_1 -1}$ satisfies this condition, we find the directional coordinates

\begin{equation}
\begin{split}
	\s_{ij} = & 
		\begin{cases}
			1, & \text{if } i = 1, j = \nu_1 \\
			-1, & \text{if } i = 1, j = \nu_1-1 \\
			0 & \text{otherwise}
		\end{cases} \\
	& \text{and} \\
	\alt{\ss} = & (0,...,0,\alt{\s}_{\l 0},...,\alt{\s}_{\l \alt{\nu}_\l})
\end{split}
\end{equation}

\noindent in $\simp{}{m}$ and $\asimp{}{m}$ respectively.

On $\simp{}{m}$ and $\asimp{}{m}$ we define polynomials of total degree $\d = |\dt|$. They are smoothly joined by enforcing the conditions in \eqref{eq:prob_cont}. With the above directional coordinates, this results in 

\begin{equation}
\label{eq:cont_simp}
\left( \c_{(\k_{11},...,\k_{1, \nu_1 -1},0)}^{(\ord)} (\ss_1) \right)_{(\k_{2},...,\k_{\l})} =
	\left( \alt{\c}_{(\k_{11},...,\k_{1, \nu_1 -1},\k_{2},...,\k_{\l-1})} \right)_{(\k_{\l 0},...,\k_{\l \nu_\l},0)}^{(\ord)} (\alt{\ss}_{\l}).
\end{equation}

To employ these conditions for the tensor-product polynomials, we first observe from lemma \ref{lm:poly_bnet} that $\bnett = \{ \ct_{\kt} \in \bnet : |\kt_i| = \dt_i, \forall i \in \{ 1,...,\l \} \}$. Reasoning from $\splo{\nu}$ we only need those conditions in \eqref{eq:cont_simp} that relate to the coefficients in $\bnett$. Applying this condition, we observe that due to $\nu_1 \neq \nu_\l$ the degrees $\dt$ get redistributed over the parts of the tensor-product polynomial. For a continuity order $\ord$, the degree in $\simp{1}{}$ reduces to $\dt_1 - \ord$ due to the elimination of $\kt_{1,\nu_1} = \ord$. At the same time this $\kt$ entry is appended to the subset $\kt_\l$, therefore increasing the degree in $\simp{\l}{}$ to $\dt_\l + \ord$.

The coefficients in $\asplo{\alt{\nu}}$ resulting from the above procedure lie on different B-nets depending on the continuity order. Because no other B-net than the one corresponding to $\dt$ are used, the conditions need to be transformed to this B-net. This requires a degree raising operation in $\asimp{\l}{}$ and a degree reduction in $\asimp{1}{}$, both of order $\ord$. After these operations the degrees $\dt$ are again equal in both simplotopes, so that the correction factors from \eqref{eq:poly_equi} can be omitted from the derivatives, resulting in \eqref{eq:cont_splo}.
\end{proof}

\begin{figure}
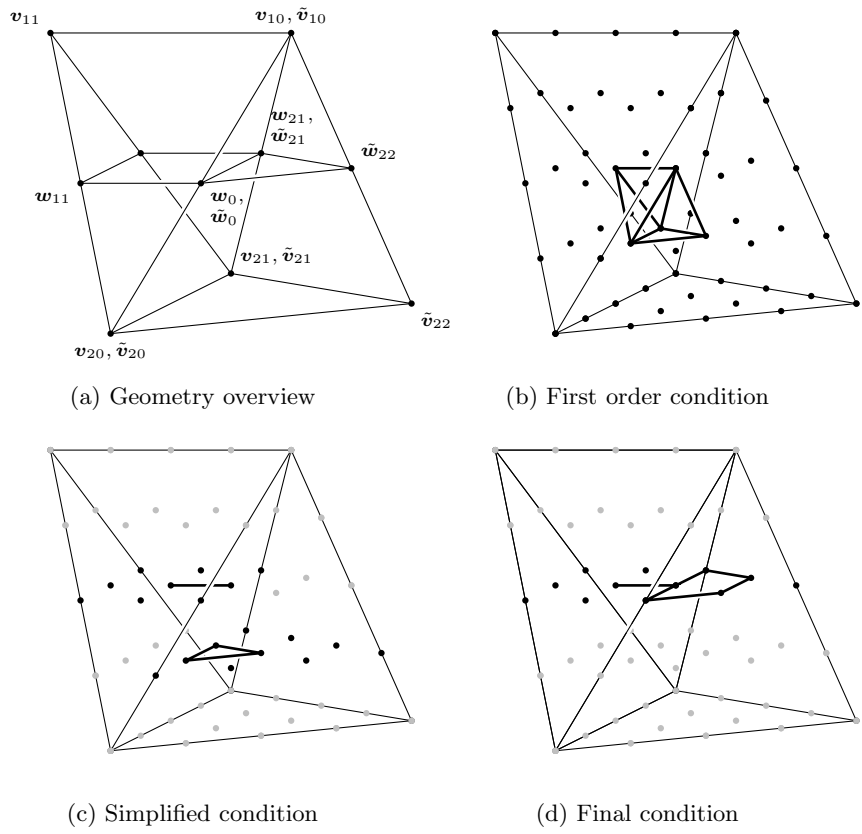

\centering
	\begin{subfigure}{0.42\textwidth}
		\begin{tikzpicture}[scale = 0.8]
			\input{Rectangle_triangle}
		\end{tikzpicture}
		\caption{Geometry overview}
		\label{fig:cont_retr_base}
	\end{subfigure} \qquad
	\begin{subfigure}{0.42\textwidth}
		\begin{tikzpicture}[scale = 0.8]
			\input{Rectangle_triangle_bnet}
		\end{tikzpicture}
		\caption{First order condition}
		\label{fig:cont_retr_bnet}
	\end{subfigure}
	\begin{subfigure}{0.42\textwidth}
		\begin{tikzpicture}[scale = 0.8]
			\input{Rectangle_triangle_splo}
		\end{tikzpicture}
		\caption{Simplified condition}
		\label{fig:cont_retr_splo}
	\end{subfigure} \qquad
	\begin{subfigure}{0.42\textwidth}
		\begin{tikzpicture}[scale = 0.8]
			\input{Rectangle_triangle_fina}
		\end{tikzpicture}
		\caption{Final condition}
		\label{fig:cont_retr_fina}
	\end{subfigure}
	\caption{Continuity conditions in the circumscribed simplices can be simplified and transformed to hold in the B-nets of the simplotopes.}
	\label{fig:cont_retr}
\end{figure}

\begin{example}[Rectangle to triangle continuity]
\label{ex:cont_retr}
Consider a rectangle $\splo{(1,1)}$ with vertices $\W = \{ \w_{10},\w_{11},\w_{20},\w_{21} \}$ and a triangle $\splo{(0,2)}$ with vertices $\alt{\W} = \{ \alt{\w}_{10},\alt{\w}_{20},\alt{\w}_{21},\alt{\w}_{22} \}$ that share an edge $\{\w_{10},\w_{20},\w_{21}\} = \{ \alt{\w}_{10},\alt{\w}_{20}, \alt{\w}_{21} \}$. As depicted in Figure \ref{fig:cont_retr_base}, we define the circumscribed simplices $\simp{}{3}$ and $\asimp{}{3}$ such that they share an edge. Note that we use $\w_0 = \w_{10} = \w_{20}$ in both simplotopes to increase the readability of Figure \ref{fig:cont_retr_base}.

We define a quartic polynomial on both simplices, and plot their B-net in Figure \ref{fig:cont_retr_bnet}. By choosing $\u = \v_{11} - \v_{10}$, the first order continuity condition depicted in the figure becomes

\begin{equation*}
\c_{(1111)} - \c_{(2011)} = 
	\alt{\s}_{20} \alt{\c}_{(1210)} + 
	\alt{\s}_{21} \alt{\c}_{(1120)} + 
	\alt{\s}_{22} \alt{\c}_{(1111)}.
\end{equation*}

\noindent The simplified condition therefore refers to only five coefficients, as displayed in Figure \ref{fig:cont_retr_splo}. Clearly we need the two B-nets highlighted in black to describe the first order continuity conditions in the combination of the rectangle and the triangle. These B-nets belong to the polynomials of degree $\dt = (2,2)$ and $\alt{\dt} = (3,1)$ respectively.

To make the condition hold for the tensor-product polynomial of degrees $\dt = (2,2)$, we need degree raising in $\asimp{2}{}$. Applying \eqref{eq:cont_rais} we find

\begin{equation*}
\begin{split}
	\c_{(1111)} - \c_{(2011)} 
	= & \frac{1}{3} \left(
		\left( \alt{\s}_{20} + \alt{\s}_{21} \right) \alt{\c}_{(1110)} + 
		\alt{\s}_{20} \alt{\c}_{(1200)} + 
		\alt{\s}_{21} \alt{\c}_{(1020)} + 
		\alt{\s}_{22} \left( \alt{\c}_{(1011)} + \alt{\c}_{(1101)} \right)
	\right),
	%
\end{split}
\end{equation*}

\noindent where we have used $|\ss| = 0$ to simplify the expression. This condition is plotted in Figure \ref{fig:cont_retr_fina}. It is clear that, due to the degree raising operation, the condition is spread out over a trapezoid in $\alt{\bnett}$.
\end{example}

The continuity conditions in \eqref{eq:cont_splo} could have also been derived without using the circumscribed simplex. Indeed, in the two- and three-dimensional cases, Lai has obtained the same results using a direct approach of setting derivatives equal on both ends of the shared facet \cite{LaiPhD}. The proposed method however, is general for all pairs of out-of-facet cospatial $n$-dimensional simplotopes of different type.



\section{Discussion}
\label{sec:conc}

In this paper we have presented continuity conditions for polynomials defined on $n$-dimensional simplotopes of different type, that share a facet and are out-of-facet cospatial. We found that these conditions can be taken from a higher-dimensional multivariate simplex spline, which provides new insights into their structure. To make all continuity conditions refer to the same B-net, we chose to redistribute the degree over simplices of the simplotope with the highest dimensional out-of-facet simplex.

In all two- and three-dimensional cases however, there is no need for degree reduction. This is because in all these cases the out-of-facet simplex $\simp{1}{}$ of $\splo{\nu}$ is one-dimensional, resulting in a zero-dimensional simplex $\asimp{1}{}$ in $\asplo{\alt{\nu}}$. In the zero-dimensional simplex, there is only one B-coefficient, of which the degree can be changed without consequences. This corresponds well with the results found by Lai and Chui \cite{LaiPhD,ChuiMu}. The first problem in which degree reduction is required is in the four-dimensional case of $\splo{(2,2)}$ and $\asplo{(1,3)}$.

The change of degree generally reduces approximation power. The optimal approximation order in the two- and three-dimensional cases has been discussed in great detail by Lai \cite{Lai}. Based on the trend in his results, it is expected that even higher degrees are required to achieve optimal approximation order in higher-dimensional problems.

It is expected that the proposed method can also be employed when the simplotopes are of equal type (that is $\nu$ = $\alt{\nu}$), and when the simplotopes are not out-of-facet cospatial. Therefore this will be our primary future research objective. To confirm our presumptions regarding the optimal approximation order in higher-dimensional problems, we are planning a detailed study of the previously described case: $\nu = (2,2)$ and $\alt{\nu} = (1,3)$.

\section{Acknowledgments}

The authors would like to thank the X-FRIS initiative of the faculty of Aerospace Engineering of the Delft University of Technology for funding this research.

\section{References}
\bibliographystyle{plain}
\bibliography{Biblio}

\end{document}